\documentclass[12pt]{amsart}

\newtheorem{theorem}{Theorem}[section]
\newtheorem{lemma}[theorem]{Lemma}

\newtheorem{proposition}[theorem]{Proposition}

\theoremstyle{definition}

\newtheorem{example}[theorem]{Example}

\newtheorem{conjecture}[theorem]{Conjecture}

\theoremstyle{remark}

\numberwithin{equation}{section}

\newcommand{\Hh}{\mathbb{H}}
\newcommand{\D}{\mathbb{D}}
\newcommand{\C}{\mathbb{C}}
\newcommand{\N}{\mathbb{N}}

\newcommand{\R}{\mathbb{R}}

\newcommand{\si}{\sigma }

\newcommand{\ga}{\gamma }
\newcommand{\Ga}{\Gamma }

\newcommand{\om}{\omega }
\newcommand{\Om}{\Omega }

\newcommand{\rea}{\operatorname{Re}}
\newcommand{\ima}{\operatorname{Im}}
\newcommand{\Arg}{\operatorname{Arg}}

\newcommand{\ntlim}{\operatorname{n.t.-lim}}

\newcommand{\Cap}{\operatorname{Cap}}
\newcommand{\bd}[1]{\partial #1}

\begin{document}
\baselineskip=18pt

\title[Pointwise convergence]{Pointwise convergence on the boundary in
the Denjoy-Wolff Theorem}

\author{Pietro Poggi-Corradini}
\address{Department of Mathematics, Cardwell Hall, Kansas State University,
Manhattan, KS 66506, USA.} \email{pietro@math.ksu.edu}
\subjclass[2000]{30D05,30C85}
\date{\today}

\begin{abstract}
If $\phi$ is an analytic selfmap of the disk (not an elliptic
automorphism) the Denjoy-Wolff Theorem predicts the existence of a point $p$
with $|p|\leq 1$ such that the iterates $\phi_{n}$ converge to $p$
uniformly on compact 
subsets of the disk. Since these iterates are bounded analytic functions,
there is a subset of the unit circle of full linear measure where they
all well-defined. We address the question of whether convergence to $p$
still holds almost everywhere on the unit circle. The answer depends on the
location of $p$ and the dynamical properties of $\phi $. We show that
when $|p|<1$(elliptic case), pointwise a.e. convergence holds if and
only if $\phi$  is not an
inner function. When $|p|=1$ things are more delicate. We show that when $\phi$
is hyperbolic or type I parabolic, then pointwise a.e. convergence holds
always. The last case, type II parabolic remains open at this moment, but
we conjecture the answer to be as in the elliptic case. 
\end{abstract}

\maketitle

\section{Introduction}\label{sec:intro}
Let $\phi$ be an analytic map defined on the unit disk $\D=\{z\in
\C:|z|<1 \}$, and assume that $\phi (\D)\subset\D$ (we call $\phi$ a
self-map of the disk from now on). The iterates of $\phi$ are
$\phi_n=\phi\circ \cdots \circ \phi$, $n$ times. The following result
is classical (elliptic automorphisms are those that can be conjugated to
a rotation).
\begin{theorem}[Denjoy-Wolff]
If a self-map of the disk $\phi$ is not an elliptic automorphism, then
there exist a point $p\in \overline{\D}$ such that the sequence
$\phi_{n} (z)$ converges uniformly on compact subsets of $\D$ to $p.$

Moreover, when $p\in \D$, $\phi (p)=p$ and $|\phi^{\prime} (p)|<1$, while
when $p\in \bd\D$, then $\phi (p)=p$ and $0<\phi^{\prime}
(p)\leq 1$ in the sense of non-tangential limits.
\end{theorem}
The point $p$ is referred to as the {\em Denjoy-Wolff point} of
$\phi$.
When $p\in \D$, the map $\phi$
is called {\em elliptic}.
When $p\in \bd\D$, $\phi$ is
called {\em hyperbolic} if $\phi^{\prime} (p)<1$
and {\em parabolic} if $\phi^{\prime} (p) =1$.

Since the functions $\phi_{n}$ are bounded analytic functions, it is
well-known that one can define the corresponding boundary functions
\[
\phi_{n}^{\star} (e^{i\theta})=\lim_{r\uparrow 1}\phi_n (re^{i\theta})
\]
for almost every $e^{i\theta}$ on $\bd\D$. More precisely, for every
$n=1,2,3,\dots$, there is a set $E_{n}\subset\D$ of linear measure
zero, so that $\phi_{n}^{\star}$ is well-defined on $\bd D\setminus E_{n}$.
Then $W=\bd \D\setminus \cup_{n=1}^{\infty}E_{n}$ has full measure and
every iterate $\phi_{n}$ extends to $W$. 

It is natural to ask whether
for almost every point $\zeta \in W$ the sequence $\phi_{n}^{\star}
(\zeta)$ still converges to $p$ as $n\rightarrow \infty$.
An answer to this question in the elliptic case (when $p$ is not
super-attracting) 
can be extracted from the proofs of our paper \cite{mich}.
The purpose of this note is
to do a more systematic study of this problem.

Before stating our result, we need to recall some definitions. 

A bounded analytic function $f$ on the unit disk $\D$ is an {\em inner
function} if the corresponding boundary function $f^{\star} (e^{i\theta})$
has modulus equal to $1$ for almost every $e^{i\theta}$ on $\bd\D$.

Recall the hyperbolic distance, defined for  $z,w\in\D$ as:
\[
\rho (z,w)=\log \frac{1+\left|
\frac{z-w}{1-\overline{w}z}\right|}{1-\left|\frac{z-w}{1-\overline{w}z}
\right|}.
\]
Given a self-map $\phi$ of parabolic type, pick a point
$z_{0}\in \D$ and let $z_{n}=\phi_n (z_{0})$ be the corresponding orbit.
Define $s_{n}:=\rho (z_{n},z_{n+1})$, i.e. the hyperbolic step of the orbit.
By Schwarz's Lemma, $\phi$ is a contraction with respect to the metric
$\rho$, in particular, $s_{n}$ is a non-increasing sequence, and hence
$s_{\infty}=\lim_{n\rightarrow \infty}s_{n}$ exists. There are two cases:
$\phi$ is {\em type I parabolic} (or {\em non-zero-step}) if
$s_{\infty}>0$, and {\em type II 
parabolic} (or {\em zero-step}) if $s_{\infty}=0$.
It follows from the main theorem of \cite{po79} that this classification
does not depend on the choice of $z_{0}$ (also see comments after
Theorem 1.8 of \cite{iber}).

\begin{theorem}\label{thm:main} 
Suppose $\phi$ is an analytic self-map
of the disk which is not an elliptic automorphism.
\begin{enumerate}
\item If $\phi$ is elliptic, then $\phi_{n}^{\star} (\zeta)$
converges to the Denjoy-Wolff point of $\phi$, for almost
every $\zeta$ in $\bd \D$, if and only if $\phi$ is
not an inner function.
\item If $\phi$ is hyperbolic or type I parabolic, then
$\phi_{n}^{\star} (\zeta)$ 
converges to the Denjoy-Wolff point of $\phi$, for almost
every $\zeta$ in $\bd \D$.
\end{enumerate}
\end{theorem}

Our proofs do not extend to the type II parabolic case, which in some
sense is more similar to the elliptic case since the hyperbolic steps
are tending to zero. In fact, we make the following conjecture.
\begin{conjecture}\label{conj}
Suppose $\phi$ is an analytic self-map
of the disk which is type II
parabolic. Then $\phi_{n}^{\star} (\zeta)$
converges to the Denjoy-Wolff point of $\phi$, for almost
every $\zeta$ in $\bd \D$, if and only if $\phi$ is
not an inner function.
\end{conjecture}
In Section \ref{sec:ell} we
tackle the elliptic case. In Section \ref{sec:bdw} we deal with the hyperbolic
and type I parabolic cases.

\section{The elliptic case}\label{sec:ell}
Assume that $\phi$ is elliptic, i.e., the Denjoy-Wolff point $p$ of
$\phi$ is in 
$\D$ (and assume that $\phi$ is not an automorphism). If $\phi$ is
inner then for almost every $\zeta\in \bd\D$, $|\phi_{n}^{\star}
(\zeta)|=1$ for all $n$, so the convergence to $p$ does not occur.

The converse is less straightforward. Assume then also that $\phi$ is
not an inner function. The argument revolves around showing that there
cannot exist a set $E\subset\bd\D$ of positive linear measure on which
all the iterates $\phi_{n}^{\star}$ have modulus one.

\subsection{An exhaustion of the unit disk}\label{ssec:exhaust}
Given $n=1,2,3,\dots$, 
fix a parameter $t>0$, and consider the open set $U_{n} (t)=\{z\in
\D:\rho (\phi_{n}(z),p)<t \}$. Then, let $\Omega_{n} (t)$ be the
connected component of $U_{n} (t)$ which contains $p$, and let $F_{n}
(t)$ be $\bd \Omega_{n} (t)\cap \D$. Notice that $F_{n} (t)$ consists
of at most countably many piecewise analytic Jordan arcs, and either there is
only one closed arc, or all the arcs have the
property that their two ends
tend to $\bd\D$ (by the maximum principle).
Let $C (r)=\{z\in\D:\rho (z,p)=r \}$, then there is $r>0$ such that
$\phi (z)\neq p$, for all $z\in C (r)$. Therefore,
we can find $t_{0}>0$ small enough so
that $\Om_{1} (t_{0})$ is compactly contained in $\D$, and therefore
$F_{1} (t_{0})$ consists of one closed Jordan arc. From now on we write
$\Omega_{n}$ for $\Omega_{n} (t_{0})$, and $F_{n}$ for $F_{n} (t_{0})$.

By definition 
\begin{equation}\label{eq:incl}
\phi_{k} (\Omega_{n+k})\subset \Omega_{n}\qquad \mbox{and}\qquad
\phi_{k}(F_{n+k})\subset F_{n} 
\end{equation}
for $n,k=1,2,3,\dots$. Moreover, we also have
\begin{equation}\label{eq:subset}
\Omega_{n}\cup F_{n}\subset 
\Omega_{n+1},
\end{equation}
for $n=1,2,3,\dots$. This because whenever $\zeta \in
\Omega_{n}\cup F_{n}$, there is a path $\ga\subset\Omega_{n}\cup
F_{n}$ connecting $p$ to $\zeta$, and by the invariant form
of Schwarz's Lemma and the fact that $\phi (p)=p$,
\[
\rho (p,\phi_{n+1} (\gamma (s) ))\leq \rho (p,\phi_{n} (\gamma (s)
))<t_{0}, \qquad 0\leq s\leq 1. 
\]

\subsection{Harmonic measure}\label{ssec:harm}
If $\Omega$ is an open set and $E$ a closed set, then we write
\[
\om (z,E,\Omega)
\]
for the Perron solution of the Dirichlet problem, in the component $U$ of
$\Omega\setminus E$ containing $z$, with data $\chi_{E}$ (the
characteristic function  of $E$). Recall that this is obtained by
taking the supremum of all the values $v (z)$, when $v$ ranges among
all subharmonic functions on $U$ such that $\limsup_{z\rightarrow
\zeta}v (z)\leq \chi_{E} (\zeta )$, for all $\zeta \in \bd\Omega\cup
E$ (these functions $v$ are often referred to as
``candidates'').

Write $\om_{n}(z)=\om (z,F_{n},\Omega_{n})$. 
We will need two results about harmonic measure.
We refer to \cite{ran} for the potential theory background that is needed. 
\begin{lemma}[Schwarz-type Lemma]\label{lem:sl}
Let $E$ be a closed set in $\D$ with $\Cap (\phi^{-1} (E))>0$. Then,
\[
\om (z,\phi^{-1} (E),\D\setminus\phi^{-1} (E))\leq\om (\phi
(z),E,\D\setminus E).  
\]
\end{lemma}
\begin{proof}[Proof of Lemma \ref{lem:sl}]
This proof is similar to the proof of Lemma 3.1 in
\cite{mich}.
Let $G=\D\setminus\phi^{-1}(E)$.
Let $v$ be a candidate for the Dirichlet problem on
$G$ with data $\chi_{\phi^{-1} (E)}$, and let
$u (z)=\om (\phi(z),E,\D\setminus E)$. When $z\in G$, 
$\phi (z)\not\in E$, and hence
$v-u$ is subharmonic on $G$. Suppose now that $\zeta \in\bd G$. There
are two cases. First assume that $\zeta \in \bd \D$, i.e. $\zeta
\not\in \phi^{-1} (E)$. Then by definition of $v$ and since $u$
is positive, $\limsup_{z\rightarrow \zeta}[v (z)-u (z)]\leq 0$.
When $\zeta\in \phi^{-1} (E)$,
$\limsup_{z\rightarrow \zeta}v (z)\leq 1$.
Also, $\Cap E>0$ by Corollary 3.6.6 of \cite{ran}, and 
at nearly every $\eta \in E$ we have
$\lim_{z\rightarrow\eta}\om (z,E,\D\setminus E)=1$, by Theorem 4.2.5 of
\cite{ran} (Kellogg's Theorem), and by Theorem 4.3.4 (b) of
\cite{ran}. Therefore,
since $\phi$ is analytic, for nearly every $\zeta \in \phi^{-1} (E)$,
$\lim_{z\rightarrow\zeta}u (z) =1$.
By the extended maximum principle for subharmonic
functions, \cite{ran} 
Theorem 3.6.9. (b), $v-u\leq 0$ on $G$, and the conclusion is reached by
taking the supremum over all the candidates $v$.  
\end{proof}

The second result is a well-known ``conditional probability
estimate''.
\begin{lemma}[Conditional Probability]\label{lem:cpe}
Suppose that $\Omega$ is an open set and $E$ is a non-empty Borel subset of
$\bd\Omega$. Also suppose that $F$ is a closed subset of $\Omega$
which separates a 
point $z\in \Omega$ from $E$ in $\Omega$, i.e., if $U$ is the
connected component of $\Omega\setminus  F$ 
containing $z$, then $E\cap \bd U=\emptyset$. We have 
\[
\om (z,E,\Omega)\leq \om (z,F,\Om)\sup_{\zeta \in F}\om (\zeta ,E,\Om).
\]
\end{lemma}
\begin{proof}
With $U$ as above, $u (w):=\om (w,F,U)=\om (w,F,\Om)$ is harmonic in
$U$. Let $v$ be a subharmonic candidate for $\om (z,E,\Omega)$. Then
$f (w):= v (w)-u (w)\sup_{\zeta \in F}\om
(\zeta ,E,\Om)$ is subharmonic in $U$. First note that since $F$
separates $z$ from $E$ in $\Omega$, we must have $\Cap F>0$, see
Corollary 3.6.4 of \cite{ran}. Then, for nearly every 
$\xi\in F$, 
$\lim_{w\rightarrow \xi}u (w)=1$, see Theorem 4.3.4 of \cite{ran}. 
Moreover, since $\xi\not\in E$
and $v$ is upper-semicontinuous, $\limsup_{w\rightarrow \xi}v (w)\leq
v (\xi)\leq \om (\xi,E,\Omega)$.
Therefore $\limsup_{w\rightarrow
\xi}f (w)\leq 0$. On the other hand, 
if $\xi\in \bd U\setminus F$, then $\xi\in \bd\Om\setminus  E$;
therefore, since $u (w)\geq 0$, and since by definition 
$\limsup_{w\rightarrow \xi}v (w)\leq 0$, we again have
$\limsup_{w\rightarrow \xi} f (w)\leq 0$.
So by the extended maximum principle for subharmonic
functions, \cite{ran} 
Theorem 3.6.9 (b), $f (w)\leq 0$ for $w\in U$. The conclusion is reached
by taking the supremum over all candidates $v$.
\end{proof}

\subsection{Non-inner elliptic selfmaps}\label{ssec:noninner}
Recall the exhaustion $\Om_{n}$ defined in Section
\ref{ssec:exhaust}, and  that $\om_{n}(z):=\om (z,F_{n},\Om_{n})$.
\begin{lemma}\label{cl:exist}
Assume that the elliptic selfmap of the disk $\phi$ is not inner. 
Then, there is an integer $N>1$ large enough such
that $\om_{N}(z)<1$ for every $z\in \Omega_{N}$. In particular,
\begin{equation}\label{eq:alpha}
\alpha :=\sup_{\zeta \in F_{1}}\om_{N}(\zeta)<1.
\end{equation}
\end{lemma}
\begin{proof}[Proof of Lemma \ref{cl:exist}]
This argument is very similar to the one on p.~506 of \cite{mich}.
We reproduce it here for convenience. Since $\phi$ is not inner, there is
a set of positive measure $A\subset W$  (recall $W\subset\bd\D$ is the set of
full-measure where all
the iterates of $\phi$ are well-defined) such that $\phi^{\star}
(A)\subset \D$. 
By Lindel\"of's Theorem (\cite{ga} p.~92), 
it is well-known that the radial limits of
$\phi$ coincide with its 
non-tangential limits. Therefore, for
$\zeta \in \bd\D$, we define the non-tangential region:
\[
\Gamma(\zeta )=\{z\in\D:|\zeta -z|< 2\frac{1+|p|}{1-|p|} (1-|z|)\}
\]
(notice that $p\in \Ga (\zeta)$
for all $\zeta \in \bd\D$).

By restricting ourselves to a subset of $A$ of positive linear measure, 
we can assume that $\sup \{\phi(z):  z\in \Ga (\zeta)\}\leq
s<1$, for some $0<s<1$ .
By uniform convergence of $\phi_{n}$ on $s\overline{\D}$, there is
$N\in \N$ such that 
$\rho (p,\phi_{N}^{\star} (z ))<t_{0}$ for all $z\in \Ga (\zeta)$ and
all $\zeta \in A$. Thus the region $G=\cup_{\zeta \in A}\Ga (\zeta)$
is a Jordan domain contained in $\Omega_{N}$. The boundary of $G$ is locally
Lipschitz, so harmonic measure on $\bd G$ is absolutely continuous
with respect to
linear measure (this follows from McMillan's Sector Theorem; see
Section 6.6 of \cite{pom}). Hence $\om (z,A,G)>0$ for $z\in G$. 
By the maximum principle, then,
$\om (z,A,\Om_{N})>0$ as well, and since $\om (z,.,\Om_N)$ is a
probability measure on $\bd\Om_{N}$, we must have $\om
(z,F_{N},\Omega_{N})<1$ for $z\in \Om_{N}$ (recall $F_{N}\subset\D$). 
\end{proof}
\begin{proposition}\label{prop:claim}
Assume that the elliptic selfmap of the disk $\phi$ is not inner. With
the notations above,
\begin{equation}\label{eq:conv}
\om_{n} (p)\rightarrow 0,\qquad \mbox{as }n\rightarrow
\infty.
\end{equation}
\end{proposition}
\begin{proof}
We apply Lemma \ref{lem:sl} with $E=F_{N}$ (where
$N$ and $\alpha $ are as in Claim \ref{cl:exist}), and
$\psi_{k}:=\phi_{(N-1)k}$  
($k=1,2,3,\dots$) instead of $\phi$, to obtain
\[
\om (z,\psi_{k}^{-1} (F_{N}),\D\setminus\psi_{k}^{-1}
(F_{N}))\leq\om (\psi_{k} 
(z),F_{N},\D\setminus F_{N})
\]
For notational simplicity, let $T_{k}=F_{N+
(N-1)k}$ and $G_{k}=\Omega_{N+ (N-1)k}$. Then,
by (\ref{eq:incl}), $T_{k}\subset \psi_k^{-1} (F_{N})$ and
$\psi_{k}^{-1} (F_{N})\cap 
G_{k}=\emptyset$. 
Therefore, for $z\in G_{k}$,
\[
\om (z,\psi_{k}^{-1} (F_{N}),\D\setminus\psi_{k}^{-1}
(F_{N}))=\om (z,T_{k},G_{k}).
\] 
Taking the supremum for $z\in T_{k-1}$, which is a subset of $G_{k}$
by (\ref{eq:subset}),  and since by (\ref{eq:incl}) $\psi_{k}
(T_{k-1})\subset F_{1}$,  we obtain:
\[
\sup_{\zeta \in T_{k-1}}\om(\zeta,T_{k},G_{k}) \leq
\sup_{\zeta \in F_{1}}\om(\zeta,F_{N},\Omega_{N})=\alpha<1. 
\]
Now, we use the conditional
probability estimate of Lemma \ref{lem:cpe}, for $n>N$ 
\[
\om_{n} (p)\leq \alpha \sup_{\zeta \in F_{N}}\om
(\zeta,F_{n},\Omega_{n}) \leq \alpha. 
\]
Likewise, for $n>2N-1$,
\[
\om_{n} (p)\leq 
\alpha  \sup_{\zeta \in F_{N}}\om
(\zeta,F_{n},\Omega_{n}) \leq \alpha \sup_{\zeta \in F_{N}}\om
(\zeta,F_{2N},\Omega_{2N}) \sup_{\zeta \in F_{2N}}\om
(\zeta,F_{n},\Omega_{n}) \leq \alpha^{2}.
\]
More generally, for $n>N+k (N-1)$, 
\[
\om_{n} (p)\leq \alpha^{k}\rightarrow 0
\]
as $k$ tends to infinity. Therefore (\ref{eq:conv}) is proved. 
\end{proof}

\subsection{Proof in the elliptic case}\label{ssec:pfint}
Observe that, given a point $\zeta \in W$, if $\phi_{n}^{\star}
(\zeta)\in \D$, then $\phi_{n+k}^{\star}(\zeta)=\phi_{k}
(\phi_{n}^{\star} (\zeta))\rightarrow p$, as $k\rightarrow \infty$. 
Thus, if the sequence
$\phi_n^{\star}$ does not converge pointwise to $p$, there is a set
$A\subset W$ of positive linear measure such that for any $\zeta \in
A$, $|\phi_{n}^{\star} (\zeta)|=1$ for all $n=1,2,3,\dots$.
We claim that
\begin{equation}\label{eq:ineq}
0<\om (p,A,\D)\leq \om_{n} (p)
\end{equation}
but letting $n$ tend to infinity and using (\ref{eq:conv}) we thereby reach a
contradiction. 

To prove (\ref{eq:ineq}), we use the fact that, although $F_{n}$ may
not separate $A$ from $p$, it at least does so ``radially''.
Fix an integer $n$, and for every $\zeta \in A$, since
$|\phi_{n}^{\star} (\zeta)|=1$, we can find $0<r 
(\zeta)<1$ so that $\rho (\phi_{n} (r\zeta),p)>t_{0}$ for $r (\zeta)\leq r<1$.
In particular, the slit $S_{\zeta}=[r (\zeta)\zeta ,\zeta )$ does not intersect
$\Omega_{n}$. So,
letting $\tilde{A}=\cup_{\zeta \in A}S_{\zeta}$, we find that
\begin{equation}\label{eq:tilde}
\om (p,\tilde{A},\D\setminus\tilde{A})\leq \om_{n} (p),
\end{equation}
as one can see from Lemma \ref{lem:cpe}, for instance.

Finally, the proof of (\ref{eq:ineq}) is completed if we can show that
\begin{equation}\label{eq:hall}
\om (p,A,\D)\leq\om (p,\tilde{A},\D\setminus\tilde{A}).
\end{equation}
To see (\ref{eq:hall}),
let $v (z)$ be a subharmonic candidate for $\om (z,A,\D)$. By the
maximum principle, $v (z)\leq 1$ for all $z\in \D$. So $v$ is also a
candidate for $\tilde{A}$ in $\D\setminus\tilde{A}$, i.e.,
$v (z)\leq \om (z,\tilde{A},\D\setminus\tilde{A})$, and (\ref{eq:hall}) is
proved by taking the supremum over the $v$'s and evaluating at $z=p$.

This completes the proof of Theorem \ref{thm:main} in the case of an
interior Denjoy-Wolff point.

\subsection{Remarks}\label{ssec:rem}
We have just shown that the pointwise a.e. convergence on $\bd
\D$ of the iterates to the Denjoy-Wolff point holds whenever the
self-map is elliptic and non-inner. As mentioned above, this fact at
least in the case when the derivative of $\phi$ at $p$ is non-zero was
already contained in \cite{mich} and \cite{arkiv}. However, there the main
tool was the conjugating map $\sigma :\D: \rightarrow \C$ due to K\oe
nigs \cite{koe}, which solves the functional equation
\[
\sigma \circ \phi (z) =\phi^{\prime} (p)\sigma (z), 
\]
and the following dichotomy was proved: either $\phi$ is not inner and
then $\sigma$ has finite non-tangential limits (in $H^{p} (\D)$)
almost everywhere on $\bd\D$, or $\phi$ is inner and then the radial
maximal function of $\sigma$ is infinite a.e. on $\bd\D$.
In the case when $\phi^{\prime} (p)=0$, one would have to use a
different conjugating map due to B\"{o}ttcher \cite{bo}, which is not
even well-defined in $\D$, but the logarithm of its modulus is, see 
\cite{cg} p.~33.
 
What we have done here is purge the map $\sigma$ from
the arguments. We will see below in the hyperbolic and parabolic
cases, that conjugations will again be useful.

\section{The case when the Denjoy-Wolff point is on the
boundary}\label{sec:bdw} 

When the Denjoy-Wolff point $p$ is on $\bd\D$ it is customary to
change variables with the M\"obius transformation $i (w+z)/
(w-z)$ so that $\phi$ becomes a self-map of the upper-half plane $\Hh$
with Denjoy-Wolff point at infinity. Julia's lemma then implies that
$\phi$ can be written as 
\begin{equation}\label{eq:standard}
\phi(z)=Az+p (z)
\end{equation}
for some $A\geq 1$ and some function $p$, with $\ima p (z)\geq 0$ for
all $z\in\Hh$, such that 
\[
\ntlim_{z\rightarrow \infty}\frac{p (z)}{z}=0.
\]
In particular, the horodisks $H (t)=\{z\in \Hh: \ima z>t \}$ ($t>0$)
are mapped into themselves, and the map $\phi $ is classified as
{\em hyperbolic} if $A>1$ and {\em parabolic} if $A=1$.

The proof of Theorem \ref{thm:main} in the elliptic case hinged on the
fact that for non-inner selfmaps of the disk there cannot exist a set
$E\subset\bd\D$ of positive linear measure on which the non-tangential limits 
$\phi_{n}^{\star}$ of each iterate $\phi_n$ all have modulus
one. This, however, is quite possible for 
non-inner selfmaps of hyperbolic and parabolic type, as the following
example shows.
\begin{example}\label{ex:hyptyp}
Let $G$ be the upper half-plane minus the slits $L_{n}=\{z=x+iy: x=2^{n},
0<y\leq 2^{n} \}$  for $n=1,2,3,\dots$, and minus the rectangle
$R=\{z=x+iy: -1\leq x\leq 1, 0<y\leq 1 \}$. The domain $G$ is simply
connected, so let $\sigma$ be the Riemann map of $\Hh$ onto $G$, such
that $\ntlim_{z\rightarrow \infty} \sigma(z)=\infty$. Defining
$\phi (z):=\sigma^{-1} (2\sigma (z))$, one can check that $\phi$ is
hyperbolic, non-inner, and all its iterates have zero imaginary part
on $\sigma^{-1} (L_{1})\subset \R$. A parabolic example can be
obtained by letting $L_{n}=\{z=x+iy:x=n, 0<y\leq 1 \}$ for
$n=1,2,3,\dots$, and $R=\{z=x+iy: x\leq 0, 0<y\leq 1 \}$.  
\end{example}
Therefore the proof in the hyperbolic and parabolic cases must
necessarily be different. We begin with the hyperbolic case.

\subsection{The hyperbolic case}\label{ssec:smh}

We need the following conjugation due to Valiron; see also \cite{jyv} for a
recent exposition of this result.
\begin{theorem}[Valiron \cite{va1}]\label{thm:valiron}
Assume $\phi$ is as in (\ref{eq:standard}) with $A>1$, i.e., $\phi$ is
hyperbolic.
Then, there is an analytic map $\si$ with $\sigma
(\Hh)\subset\Hh$ such that $\ntlim_{z\rightarrow \infty }\sigma
(z)=\infty$  and actually $\sigma$ is isogonal at infinity, i.e.
$\ntlim_{z\rightarrow \infty}\Arg \sigma (z)=0$, 
which satisfies the functional equation:
\[
\sigma \circ \phi =A\sigma 
\]
\end{theorem}
Since $\sigma$ is bounded analytic (after a change of variables), it has
non-tangential limits almost everywhere on the real axis. We let
$\sigma^{\star}$ be the boundary function.
By the F. and M. Riesz Theorem (\cite{ga} p.~65) the set $Z=\{x\in \R:
\sigma^{\star} (x)=0\}$ has measure zero.
Suppose $x\in \R\setminus Z$ 
is a point at which $\sigma^{\star} (x)$ and all the iterates
$\phi_{n}^{\star} 
(x)$ are well-defined, and assume that  
\[
\liminf_{n\rightarrow \infty}|\phi_{n}^{\star} (x)|=R<\infty. 
\]
Choose a sequence of integers $n$ for which $|\phi_{n}^{\star} (x)|<2R$.
For each such $n$, define $\ga_{n} (t)=\phi_{n} (x+it)$, for $t>0$.
Then, $\ga_{n}$ is a curve which connects
the ball $B (2R)=\{|z|\leq 2R \}$ to infinity.
By the
F. and M. Riesz Theorem, we also have that 
the boundary function $\sigma^{\star} (x)$ is
finite almost everywhere. In particular, 
there exists
$s>2R$ such that $\si^{\star} (s),\sigma^{\star}
(-s)<\infty$. Therefore if $\Ga =\{z\in \Hh: |z|=s \}$, 
then 
\begin{equation}\label{eq:mr}
M:=\sup_{z\in \Ga}|\sigma(z)|<\infty.
\end{equation}
Let $z_{n}=\phi_{n} (x+it_{n})$ be a
point in $\Ga \cap \ga_{n}$, which is necessarily non-empty.
Then,
\[
M\geq|\sigma (z_{n})|=|\sigma (\phi_{n} (x+it_{n}))|=A^{n}|\si (x+it_{n})|
\]
Therefore, since $A>1$,
 $\lim_{n\rightarrow \infty}|\si (x+it_{n})|=0$. However,
since $\lim_{t\rightarrow +\infty}\sigma (x+it)=\infty$, and $|\si
(x+it)|>0$ for all $t>0$, we must conclude that $t_{n}$ tend to $0$,
i.e., that $\sigma^{\star} (x)=0$. But this contradicts our hypothesis
$x\not\in Z$.

In conclusion, we find that except for a set of linear measure zero,
at all points $x$ where the iterates $\phi_{n}^{\star} (x)$ are well-defined
we have $\lim_{n\rightarrow \infty}|\phi_{n}^{\star} (x)|=\infty$,
which proves Theorem \ref{thm:main} in the hyperbolic case.

\subsection{The type I parabolic case}\label{ssec:typeIp}

The counterpart of Valiron's Theorem in the type I parabolic case is
the following result of Pommerenke.
\begin{theorem}[Pommerenke, \cite{po79} Theorem 1 and (3.17)]\label{thm:pom}
Let $\phi $ be an analytic self-map of $\Hh $ of parabolic type as in
(\ref{eq:standard}) with $A=1$, and 
let $\{z_{n}=\phi _{n} (i)\}_{n=0}^{\infty }$ be a forward-iteration
sequence.  Then, 
\[
\frac{\ima z_{n+1}}{\ima z_{n}}\longrightarrow 1
\]
as $n$ tends to infinity.

Moreover, if $\phi$ is type I, i.e. $\rho
(z_{n},z_{n+1})\downarrow s_{\infty}>0$, 
letting $z_{n}=u_{n}+iv_{n}$ and considering the
automorphisms of $\Hh $ given by $M_{n} (z)= (z-u_{n})/v_{n}$,
the normalized iterates 
$M_{n}\circ \phi _{n}$ converge uniformly on compact subsets of $\Hh $
to a function $\sigma $ which satisfies the functional equation
\[
\sigma \circ \phi =\sigma + b
\]
where 
\begin{equation}\label{eq:b}
b:=\lim_{n\rightarrow \infty }\frac{u_{n+1}-u_{n}}{v_{n}}\neq 0.
\end{equation}
\end{theorem}
The conjugation $\sigma$ is a selfmap of $\Hh$ by construction and
Pommerenke also shows that $\ntlim_{z\rightarrow \infty}\ima \sigma
(z)=+\infty$ (and actually the region of convergence can be extended
to a tangential one). However, this is not enough information about the
behavior of $\sigma$ at infinity, in particular some information on
the behavior of $\rea \sigma (z)$ at infinity is necessary if
one wants to repeat the same argument as in the hyperbolic case.

Instead we modify the argument slightly.
Let $\sigma$ be the conjugation of Theorem \ref{thm:pom} and assume
without loss of generality that the constant $b$ in (\ref{eq:b}) is
positive, so that $u_{n}$ is eventually an increasing sequence and
since $v_{n}\geq v_{0}$: 
\begin{equation}\label{eq:rea}
\lim_{n\rightarrow
\infty}u_{n}=+\infty.
\end{equation}
Let $\sigma^{\star}$ be the boundary function. Suppose $x\in \R$ is a
point at which all the iterates $\phi_{n}^{\star} (x)$ are
well-defined, and where $\si^{\star} (x)$ is finite.
Assume also that
\[
\liminf_{n\rightarrow \infty}|\phi_{n}^{\star} (x)|=R<\infty. 
\]
Instead of considering the half-line $\{z=x+it, t>0 \}$, 
let $z_{n}=\phi_{n} (i)=u_{n}+iv_{n}$, and define $P$ to be the
polygonal curve 
\[
P=[x,i]\cup[i,z_{1}]\cup[z_{1},z_{2}]\cup\cdots
\]
At one end, $P$ tends non-tangentially to $x$.
Near infinity, $P$ is a simple curve tending to infinity. Moreover, by
(\ref{eq:rea}), 
$\rea z\rightarrow +\infty$ as $z$ tends to infinity along $P$. Also,
as Pommerenke remarks, in \cite{po79} Remark 1, 
\begin{equation}\label{eq:rem}
\frac{v_{n}}{u_{n}}\rightarrow 0\qquad \mbox{ as }n\rightarrow \infty, 
\end{equation}
so the argument of $z$ tends to zero as $z$ tends to infinity along $P$. 

Choose a sequence of integers  $n$ so that $|\phi_{n}^{\star}
(x)|<2R$, and let $\ga_{n}=\phi_{n} (P)$. By construction, if
$z$ tends to $x$ along $P$, then $\phi_{n} (z)$ tends to
$\phi_{n}^{\star} (x)$, and hence intersects the ball $B (2R)$. 
If $z$ tends to infinity
along $P$, we claim that 
\begin{equation}\label{eq:P}
\lim_{P\ni z\rightarrow \infty}|\phi_{n}(z)|=+\infty. 
\end{equation}
In fact, find $k$ so that $z\in [z_{k},z_{k+1}]$; then by Schwarz's
Lemma 
\[
\rho (\phi_{n} (z),z_{k+n})\leq \rho (z,z_{k})\leq \rho
(z_{1},z_{2})<\infty,
\]
so 
$\phi_{n} (z)$ tends to infinity.

Also, we claim that
\begin{equation}\label{eq:inf}
m:=\inf_{z\in P}\rea \sigma (z)>-\infty 
\end{equation}
This holds on $[x,i]$ by our choice of $x$. Also, for $k=1,2,3,\dots$,
$\sigma (z_{k})=\sigma (\phi_{k} (i)) =\sigma (i)+kb$, and for $z\in
[z_{k},z_{k+1}]$; so 
Schwarz's Lemma implies
\[
\rho (\sigma (z),\sigma (i)+kb)\leq  \rho (z,z_{k})\leq \rho
(z_{1},z_{2})<\infty
\]
which yields (\ref{eq:inf}) at once.

Find $s>2R$ such that $\sigma^{\star} (s),\sigma^{\star} (-s)<\infty$,
which can be done since $\sigma$ is a selfmap of $\Hh$. Then each
curve $\ga_{n}$ must intersect the circle $\{|z|=s \}$ at a point
of the form $\phi_{n} (w_{n})$ for some $w_{n}\in P$. Also,
$\sup_{|z|=s}|\sigma(z)|=M<\infty$, so by (\ref{eq:inf}) 
\[
M\geq |\sigma (\phi_{n}
(w_{n}))|=|\sigma (w_{n})+nb|\geq \rea \sigma (w_{n}) +nb\geq
m+nb\rightarrow+\infty  
\]
which is a contradiction.

Thus,
except for a set of linear measure zero,
at all points $x$ where the iterates $\phi_{n}^{\star} (x)$ are well-defined
we have $\lim_{n\rightarrow \infty}|\phi_{n}^{\star} (x)|=\infty$,
which proves Theorem \ref{thm:main} in the type I parabolic case.

\section{Remarks about the type II parabolic case}\label{sec:id}

Here we collect some remarks on the type II parabolic case, or, in
Pommerenke's terminology, the {\em identity case}. Recall that these
are analytic self-maps $\phi$ of the upper half-plane $\Hh$ that can be
written as in (\ref{eq:standard}) with $A=1$, and such that the
hyperbolic steps $\rho (z_{n},z_{n+1})$ of the forward-iteration
sequence $z_{n}=\phi_{n} (i)$ tend to zero. We have already mentioned
in the introduction that the fact that the hyperbolic steps tend to
zero does not depend on the choice of the starting point $i$. Moreover,
given any point $z\in \Hh$, we also have that $\rho (\phi_{n}
(z),\phi_n (i))\rightarrow 0$, as $n$ tends to infinity. This also
follows from Pommerenke's Theorem 1 \cite{po79}. In fact, with the
notations of Theorem \ref{thm:pom}, the normalized
iterates $M_{n}\circ \phi_{n}$ converge uniformly on compact subsets
of $\Hh$ to $i$ in this case. So 
\begin{equation}\label{eq:zero}
\rho (\phi_{n}(z),\phi_n (i))=\rho (M_{n}\circ \phi_{n}(z),M_{n}\circ
\phi_n (i))= \rho (M_{n}\circ \phi_{n}(z), i)\rightarrow 0. 
\end{equation}
It also follows from (\ref{eq:standard}) that for any given $z\in \Hh$
the sequence of imaginary parts $\ima \phi_{n} (z)$ is strictly
increasing, hence it either has a finite limit or it tends to
infinity. Again this is a property that does not depend on $z$, and in
the type I parabolic case both cases arise.
In \cite{iber}, we left open the problem of producing examples in the
type II parabolic case as well. However, we now realize that this question
can be easily answered.
\begin{proposition}\label{prop:typeii}
If $\phi$ is a type II parabolic self-map of $\Hh$ as above, then
given $z\in \Hh$,
\[
\lim_{n\rightarrow \infty}\ima \phi_{n} (z)=+\infty. 
\] 
\end{proposition}
The proof of this proposition is immediate because otherwise one would have
\[
\ell_{\infty}:=\lim_{n\rightarrow \infty}\ima \phi_{n} (i)<\infty 
\]
and by (\ref{eq:zero}) $\lim_{n\rightarrow \infty}\ima \phi_{n}
(z)=\ell_{\infty}$ as well, for any $z\in \Hh$, which contradicts the
fact that $\ima \phi_n (z)$ increases as soon as $\ima z>\ell_{\infty}$.

\bibliographystyle{alpha}

\end{document}